\documentclass[12pt]{article}
\usepackage[final]{epsfig}
\usepackage{graphics}
\usepackage{amsmath}
\usepackage{amsfonts}
\usepackage{latexsym}
\usepackage{amssymb}
\usepackage{graphicx}
\usepackage{epstopdf}
\usepackage{graphicx,epsfig,psfrag}
\newtheorem{lemma}{Lemma}

\newtheorem{theorem}[lemma]{Theorem}

\begin{document}
\newcommand{\eps}{{\varepsilon}}
\newcommand{\proofend}{$\Box$\bigskip}
\newcommand{\C}{{\mathbf C}}
\newcommand{\Q}{{\mathbf Q}}
\newcommand{\R}{{\mathbf R}}
\newcommand{\Z}{{\mathbf Z}}
\newcommand{\RP}{{\mathbf {RP}}}
\newcommand{\CP}{{\mathbf {CP}}}
\newcommand{\A}{{\rm Area}}
\newcommand{\Le}{{\rm Length}}

%%%%%%%%%%%%%%%%%%%%%%%%%%%%%%%%%%%%%%%%%%%%%%

\newcommand{\marginnote}[1]
{%\mbox{}\marginpar{\center{\hspace{0pt}\tiny{\bf#1}}}
}
\newcounter{ml}
\newcommand{\bk}[1]
{\stepcounter{ml}$^{\bf ML\thebk}$%
\footnotetext{\hspace{-3.7mm}$^{\blacksquare\!\blacksquare}$
{\bf ML\thebk:~}#1}}

\newcounter{st}
\newcommand{\st}[1]
{\stepcounter{st}$^{\bf ST\thest}$%
\footnotetext{\hspace{-3.7mm}$^{\blacksquare\!\blacksquare}$
{\bf ST\thest:~}#1}}

%%%%%%%%%%%%%%%%%%%%%%%%%%%%%%%%%%%%%%%%%%%%%%%

\title {Schr\"odinger's equation and ``bike tracks" -- a connection.}
\author{Mark Levi\thanks{
Department of Mathematics,
Pennsylvania State University, University Park, PA 16802, USA;
e-mail: \tt{levi@math.psu.edu. }
}
\\
}
 \bibliographystyle{plain}
\maketitle
 
 The purpose of this note is to demonstrate an equivalence between two classes of objects: the stationary Schr\"odinger equation on the one hand, and the ``bicycle tracks" on the other. We begin with the description of the latter. 

\begin{figure}[thb]
	\center{  \includegraphics{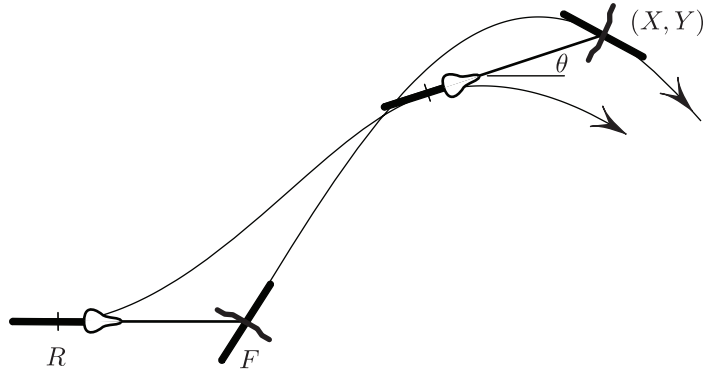}}
	\caption{An idealized bike. In this example $F$ travels along an arc of a sine curve.}
	\label{fig:bike}
\end{figure}

 A (very) idealized model of a bicycle, shown in Figure~\ref{fig:bike}, is a segment $ RF$ of constant length which is allowed to move in the plane as follows: the path of the ``front" $F$ is prescribed, while the   velocity of the ``rear" $R$ is constrained to the line 
 $ RF $: the ``rear wheel" does not sideslip.  
  %%%?  
%%%  
 If $ (X(t),Y(t)) $ is a parametric representation of the motion of $F$  then the angle  $\theta$  between $ RF $ and the $x$--axis in the plane satisfies the differential equation 
 \begin{equation} 
	\dot \theta = \dot Y  \cos \theta - \dot X \sin \theta , 
	\label{eq:a}
\end{equation}   
expressing the fact that infinitesimal displacement of $R$ is aligned with the direction $ e^{i \theta} $ of the segment. 
  
\begin{figure}[thb]
	\center{  \includegraphics{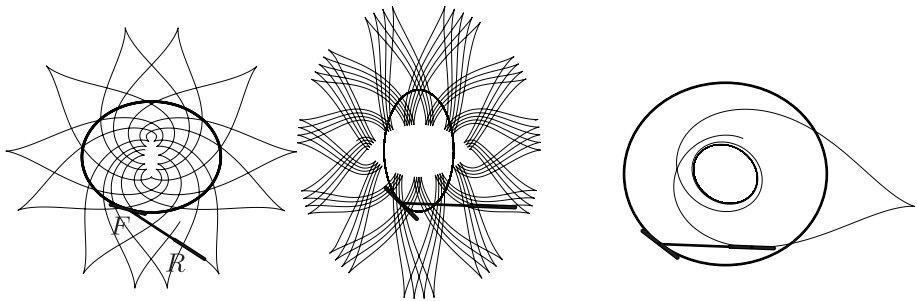}}
	\caption{The path traced out by the rear wheel as the front wheel  repeatedly traverses  a closed curve.}
	\label{fig:somebiketracks}
\end{figure}

 Some examples of tracks are given in Figure~\ref{fig:somebiketracks}. 
 \vskip 0.1 in 
%%%%%      
 \paragraph {A very brief history.} The  idealized ``bike" of Figure~\ref{fig:bike} has been studied since the second half of 19th century (see \cite{foote} and references therein), and up  
 to the present time (\cite{foote, Levi-Tabachnikov}). It was observed that the ``bike"  arises as an asymptotic limit  of a system describing a particle in  a rapidly oscillating potential; it is interesting that  the nonholonomic ``bike" is   a singular limit of a holonomic system (the details can be found in  \cite{levi-weckesser}, and in \cite{levi2}). 

\paragraph {Stationary Schr\"odinger's equation} 
 \begin{equation} 
	 \ddot x + p(t) x = 0
	\label{eq:s}
\end{equation}    
	 is a classical object of mathematical physics,  arising  in many   settings  in mathematics,   physics and engineering. This system  has been studied for nearly two centuries.  Known also as  Hill's equation, it comes up in studying
  the spectrum of hydrogen atom, in celestial mechanics \cite{siegel-moser},  in particle accelerators \cite{wiedemann}, in forced vibrations,  in wave propagation, and in many more problems. Hill's operator deforms isospectrally when its potential evolves under  the Korteweg--de Vries(KdV) equation, thus providing an explanation of complete integrability of the latter \cite{kruskal}, \cite{lax}, \cite{novikov}.  
  The  1989 Nobel Prize in physics was awarded to W. Paul for his invention of  an electromagnetic trap, now called the Paul trap,  used to suspend charged particles. The mathematical substance of Paul's discovery amounts to an observation on  Hill's equation, as expained in Paul's Nobel lecture \cite{paul1}. Incidentally, \cite{levi} contains a geometrical explanation, as an alternative to Paul's analytical one,   of why the trap works.   Stability of the famous Kapitsa pendulum \cite{arnold, kapitsa}  
is also explained by the properties of Hill's equation (Stephenson gave an experimental demonstration of stability of the so--called  Kapitsa pendulum in 1908 \cite{stephenson}, about half a century before Kapitsa's paper).  
The long history of Hill's equation is reflected in the rich body of classical literature of the 18th and 19th centuries on the eigenfunctions of special second order equations (polynomials of Lagrange, Laguerre, Chebyshev, Airy's function, etc.), to the more recent work on inverse scattering and on geometry of  ``Arnold tongues"  \cite{PS, gelfand-levitan, arnold7, levy-keller, WK1, WK2, marchenko, novikov, broer-levi, broer-simo}, \cite{levy-keller},  \cite{arnold7}, \cite{WK1,WK2}.

 \vskip 0.3 in

  \section{The main result}
\begin{figure}[thb]
	\center{  \includegraphics{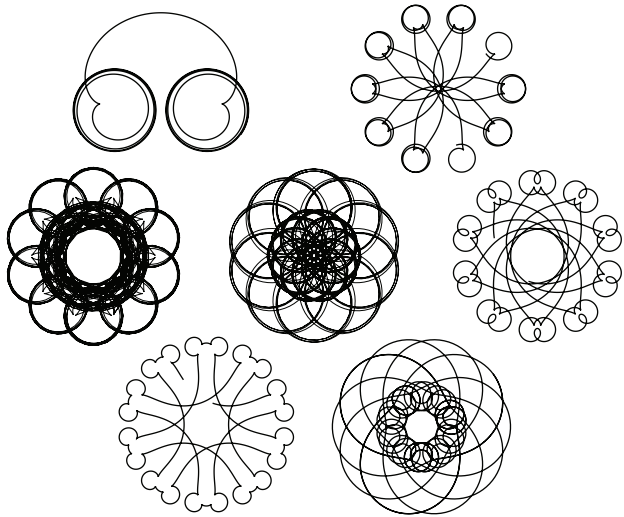}}
	\caption{Front wheel paths representing different potentials. The top left corresponds to the solitary wave of the KdV. The remaining ones are generated by various trigonometric potentials. Alternatively, these paths are trajectories of the particle subject to the ``magnetic" force   (\ref{eq:vp})  with different prescribed speeds $ v=v(t) $. }
	\label{fig:variouspotentials}
\end{figure}

\begin{theorem} \label{thm:main} Let a Schr\"odinger potential $ p(t) $ in   (\ref{eq:s}) be given. We associate with  $ p $ the   front wheel path $ (X(t), Y(t)) $ as follows:  defining
 \begin{equation} 
	   \varphi (t) = t+ \int_{0}^{t} p( s ) \,ds, 
	\label{eq:vphi}
\end{equation}   
we set 
\begin{equation}  
   \left\{ \begin{array}{l} 
    X(t) =  - \int_{0}^{t}(1-p(\tau)) \sin \varphi (\tau) \;d\tau  \\[3pt] 
    Y(t) = \  \ \    \int_{0}^{t}(1-p(\tau)) \cos  \varphi (\tau)  \;d\tau .\end{array} \right.     
	\label{eq:fwp}
\end{equation}  
 If the potential $ p $ and the path $ (X,Y) $ are thus related, then the two problems: the corresponding Schr\"odinger equation   (\ref{eq:s}) and the bike problem     (\ref{eq:a})     are equivalent in the sense that  
   \begin{equation} 
	\theta = 2 \arg (x+i \dot x ) +\varphi, 
	\label{eq:relation}
\end{equation}   
  where $ \varphi $ is given by   (\ref{eq:vphi}).\footnote{More precisely, if  (\ref{eq:relation}) holds for $ t=0 $, then it holds for all $t$.}
 \end{theorem} 
 Figure~\ref{fig:variouspotentials} shows paths corresponding to various potentials. 
 %%%?  
%%%  
 
 \subsection{A  reformulation of the main result.} 
 The track    (\ref{eq:fwp}) can be thought of as the path of a particle subject to a strange magnetic--like force defined in the next paragraph.  
 
\begin{figure}[thb]
	\center{  \includegraphics{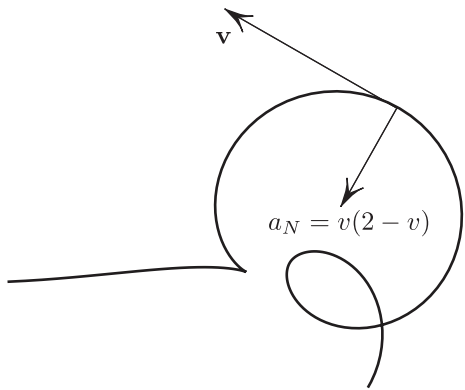}}
	\caption{Trajectory of a particle subject to force (\ref{eq:magforce}) with $v$ prescribed. At the cusp $v$ changes sign.  }
	\label{fig:magforce}
\end{figure}

  \paragraph {A pseudo--magnetic force} 
 Let $v=v(t) $ be a given function of time, and   consider a point mass $ m=1 $  moving in the plane with  speed  $v$  and subject to normal acceleration due to the following magnetic--like force: 
  \begin{equation} 
	{\bf F} ={\bf a}_\perp=  i\,(2-v){\bf v} 
	\label{eq:magforce}
\end{equation}     
acting normal to the velocity ${\bf v}$. Note that the tangential velocity $v$ is prescribed (one can imagine a tangential force acting on the particle in addition to the normal force (\ref{eq:magforce})), and that  the normal acceleration is slaved to $v$.
We allow $v$ to change sign, so that $ v = \pm | {\bf v} | $; if $v$ changes sign,   the particle reverses the direction of motion, as illustrated in   Figure~\ref{fig:magforce}.  %%%%%  

\vskip 0.1 in   
The main result can now be reformulated as follows. 

\begin{theorem} \label{thm:magnetic}
Consider the Schr\"odinger equation   (\ref{eq:s}) with potential $p(t)$. Define 
\begin{equation} 
	v(t)=  1-p(t), 
	\label{eq:vp}
\end{equation} 
and let   $ (X(t), Y(t)) $ be a path of the ``magnetic" particle defined in the preceding paragraph. Then the Schr\"odinger equation   (\ref{eq:s}) and the bike problem   (\ref{eq:a}) are equivalent, i.e. they transform into one another via the transformation  
\begin{equation} 
	\theta = 2 \arg (x+i \dot x ) + \varphi, 
	\label{eq:relation1}
\end{equation}       
	where $ \varphi = t+ \int_{0}^{t} p( \tau ) \,d\tau $. 
\end{theorem}

  \vskip 3 in 
  \section{Proofs. } 
  \noindent{\bf Proof of Theorem \ref{thm:main}.} 
We begin by writing the  Schr\"odiner  equation   (\ref{eq:s}) as a system 
\begin{equation}  
   \left\{ \begin{array}{l} 
\dot x = y  \\[3pt] 
\dot y  = - p(t) x \end{array} \right.   
\end{equation}  
or in matrix form 
\begin{equation} 
	 \dot z = P(t) z, \  \  \hbox{with}  \  \  P = \left( \begin{array}{cc} 0 & 1 \\ -p & 0 \end{array} \right) .
	 \label{eq:hillsvector}
\end{equation}   
The main point of the proof is to observe that Schr\"odiner  equation   (\ref{eq:hillsvector}) in a rotating frame becomes equivalent to the Ricatti equation for the bicycle. To make this precise, let   
\begin{equation} 
	\psi = \psi (t) = - \frac{1}{2}\biggl( t+ \int_{0}^{t} p( \tau ) d\tau \biggl);   
	\label{eq:psi}
\end{equation}   
note that $\dot x \psi $ is half the curl of the vector field in   (\ref{eq:hillsvector}), i.e.  the average angular velocity of the vector field around the origin.  
Introduce the rotation through angle $\psi$: 
\begin{equation} 
	R= R(t)    = \left( \begin{array}{cc} \cos \psi  & - \sin \psi  \\ \sin \psi  & \  \ \   \cos \psi  \end{array} \right). 	\label{eq:rotation}
\end{equation}   
  
To rewrite the Schr\"odinger equation  (\ref{eq:hillsvector}) in the rotating frame we  introduce the new variable $w$ via  
\begin{equation} 
	z= R_\psi w. 
	\label{eq:zwchange}
\end{equation}   
We obtain a new system equivalent to   (\ref{eq:hillsvector}): 
 \begin{equation} 
	\dot w = (R ^{-1} P R - R ^{-1} \dot R) w .  
	\label{eq:w}	
\end{equation} 
A computation confirms the expectation that  coefficient matrix of this system should be symmetric (since we cancelled angular velocity) and  traceless (since the transformation is area--preserving and since the original matrix was traceless):  
\begin{equation} 
	 R ^{-1} P R - R ^{-1} \dot R = 
	\left( \begin{array}{cc} r & \  \  s \\ s & -r \end{array} \right),  
	\label{eq:mw}
\end{equation} 
where  
\begin{equation}  
  	   \begin{array}{l} 
  	  r= \frac{1}{2} (1-p)\sin  2 \psi   \\[3pt] 
 	   s= \frac{1}{2} (1-p)\cos 2 \psi. 
 	\end{array}     
		 \label{eq:rs}
\end{equation}  
According to   (\ref{eq:zwchange}), we have 
\begin{equation} 
	\arg (x+iy) = \arg w + \psi, 
	\label{eq:zwangles}
\end{equation}   
and we now show that $ \theta = 2 \arg w $ satisfies the bicycle equation   (\ref{eq:a}); this would complete the proof. Indeed, then     (\ref{eq:zwangles})  would become 
\[
	\arg (x+iy) = \frac{1}{2}  \theta  + \psi, \  \  \hbox{or}  \  \  \theta = 2 \arg z- 2 \psi,  
\] 
which indeed coincides with   (\ref{eq:relation}) since $-2 \psi = \varphi$ according to   (\ref{eq:psi})  and   (\ref{eq:vphi}). 

 To derive the equation for $ \arg w $ the  we write our  system    (\ref{eq:w})-(\ref{eq:mw})  for $w$ explicitly: 
\begin{equation}  
   \left\{ \begin{array}{l} 
   \dot u = ru+sv \\[3pt] 
   \dot v=su-rv,  \end{array} \right.   
   \label{eq:gl}
\end{equation}  
and let $\widehat w =  \arg w =  \arg (u+iv) $.\footnote{the wedge in 
$\widehat w  $ reminds of the angle.} Now
\[
	\frac{d}{dt} \widehat w = \frac{\dot v u - \dot u v}{u ^2 + v ^2 } 
	\buildrel{  (\ref{eq:gl})}\over{=}  \frac{s u ^2 -2r uv-s v ^2 }{u ^2 + v ^2},  
\] 
so that 
\[
	\frac{d}{dt} \widehat w = s   \cos ^2 \widehat w - 2r \cos   \widehat w \sin  \widehat w - s \sin ^2 \widehat w
\] 
This can be rewritten in terms of double angle $ 2 \widehat w $ as follows: 
\begin{equation} 
	  \frac{d}{dt}(2\widehat w) = 2s \cos 2 \widehat w - 2r \sin 2 \widehat w.  
	\label{eq:2alpha}
\end{equation}   
This ODE is identical to the the bicycle equation:   
\[
	\dot \theta = \dot Y \cos \theta - \dot X \sin \theta  
\] 
provided we set 
\[
	\dot X=  2r , \  \  \dot Y =  2s, 
\] 
or, recalling  the definition  (\ref{eq:rs}) of $r$ and $s$, provided
  
\begin{equation}  
     \begin{array}{l} 
    \dot X = (1-p)\sin  2 \psi=- (1-p)\sin \varphi ,   \\[3pt] 
    \dot Y = (1-p)\cos 2 \psi = \  \  (1-p)\cos \varphi.  \end{array}  
    \label{eq:XY}
\end{equation}    
We conclude that $ \theta = 2 \widehat w $ in the sense that both angles satisfy the same differential equation provided we define $ X, \  Y $ by   (\ref{eq:XY}) or   (\ref{eq:fwp}). This completes the proof. 
 \hfill $\diamondsuit$  
\vskip 0.3 in 
\noindent{\bf Proof of Theorem \ref{thm:magnetic}.}
Consider  the motion $ (X(t), Y(t)) $ given by  (\ref{eq:fwp}). The velocity of this motion
\begin{equation}  
   \left\{ \begin{array}{l} 
    \dot X = -(1-p(t)) \sin \varphi   \\[3pt] 
    \dot Y = \ \ (1-p(t))\cos \varphi  \end{array} \right.   
\end{equation}  
The speed   $ v=1-p$ is in the direction $ \varphi + \pi /2 $ if $ 1-p>0 $ and in the opposite direction if $ 1-p<0 $. The angular velocity of this motion is 
\[
	\omega = \dot \varphi = 1+ p, 
\] 
and thus the normal acceleration 
\[
	a_\perp = \omega v . 
\] 
But 
\[
	\omega =1+p   \buildrel{  (\ref{eq:vp})}\over{=} 1+(1-v)=2-v,
\] 
so that $ a_\perp = v(2-v)$. 
 \hfill $\diamondsuit$  
  \bibliography{master}
 
 {\bf Acknowledgments}.  The author's research was supported by an NSF grant   DMS-0605878. 
 \end{document}